\documentclass{article}
\usepackage{amsfonts,amssymb,mathrsfs,amsmath}
\usepackage[english]{babel}
\usepackage[dvips]{color}
\usepackage{graphicx}
\title{\Large{The Coin Exchange Problem and the Structure of Cube Tilings}}
\date{}
\author{Andrzej P. Kisielewicz and Krzysztof Przes{\l}awski\\
\\
{\small Wydzia{\l} Matematyki, Informatyki i Ekonometrii, Uniwersytet Zielonog\'orski}\\
{\small ul. Z. Szafrana 4a, 65-516 Zielona G\'ora, Poland}\\
{\small A.Kisielewicz@wmie.uz.zgora.pl}\\
{\small K.Przeslawski@wmie.uz.zgora.pl}}
\newtheorem{pr}{\sc Proposition}

\newtheorem{tw}[pr]{\sc Theorem}

\newtheorem{df}{\sc Definition}

\newtheorem{uw}{\sc Remark}
\newtheorem{uwi}[uw]{\sc Remarks}

\newtheorem{nap}{\sc Example }

\newtheorem{nps}[nap]{\sc Examples}

\def\kal #1 #2{\mathscr{#1}^{#2}}

\def\skok{\vskip 0.1in}

\def\zet{\mathbb{Z}}

\def\er{\mathbb{R}}

\def\te{\mathbb{T}}

\def\przyst #1{\operatorname{mod} #1}

\def\gruba #1{\boldsymbol{#1}}
\begin{document}
\maketitle
\begin{abstract} 
Let $k_1,\ldots, k_d$ be positive integers and let $D$ be a subset of $[k_1]\times\cdots\times [k_d]$, whose complement can be decomposed into disjoint sets of the form $\{x_1\}\times\cdots\times \{x_{s-1}\}\times [k_s]\times \{x_{s+1}\}\times\cdots\times \{x_d\}$. We conjecture that the number of elements of $D$ can be represented as a linear combination of the numbers $k_1,\ldots, k_d$ with non-negative integer coefficients. A connexion of this conjecture with the structure of periodical cube tilings is revealed. 

\medskip
\noindent
\end{abstract}
For any positive integer $n$, we denote by $[n]$ the set $\{1,\ldots, n\}$. We extend this notation to vectors $\gruba k=(k_1,\ldots, k_d)$ with positive integer coordinates: $[\gruba k]:=[k_1]\times\cdots\times[k_d]$. If all $k_i$ are greater than 1, then $[\gruba k]$ is said to be a \textit{(discrete) $d$-box}. A \textit{line} in $[\gruba k]$ is any set of the form
$$
\{x_1\}\times\cdots\times \{x_{s-1}\}\times [k_s]\times \{x_{s+1}\}\times\cdots\times \{x_d\},    
$$
where $s\in [d]$, and $x_i\in [k_i]$. A subset $D$ of $[\gruba k]$ is said to be \textit{complementable by lines} if its complement $[\gruba k]\setminus D$ can be represented as a union of disjoint lines.  

A non-negative integer $n$ is \textit{representable} by $\gruba k$ if there are non-negative integers $n_1,\ldots, n_d$ such that 
$$
n=n_1k_1+\cdots+n_dk_d.    
$$
In other words, the amount $n$ can be changed using coins of denominations $k_1,\ldots, k_d$. As a consequence of this interpretation, the problem of representability is often called \textit{the coin exchange problem}.  

 The following conjecture arises from certain problems concerning periodical cube tilings, as we shall explain it  later on.
\skok
\noindent
\textbf{Conjecture} 
\begin{itemize}
\item[]
For each  $d$-box $[\gruba k]$, if $D\subseteq [\gruba k]$ is comlementable by lines, then the size $|D|$ of $D$ is representable by $\gruba k$. 
\end{itemize}
It is not difficult to confirm this conjecture  for $\gruba k=(m,\ldots,m,n)$, where $m$ and $n$ are arbitrary positive integers. 
If $d=3$, then verification of the conjecture reduces to a strictly numerical problem:
\begin{itemize}
\item[]
Show that for every positive integers $1<k_1<k_2<k_3$, and $1\le l_i\le k_i-1$, $i=1,2,3$,  the number   
$$
l_1l_2l_3+(k_1-l_1)(k_2-l_2)(k_3-l_3)
$$
is representable by $(k_1,k_2,k_3)$.
\end{itemize}
  
This problem has been tested for a wide range of data by M. Ha\l uszczak and, independently, by A. Zieli\'nski, an MSc student of the second author. In particular, it has been tested for all $1<k_1<k_2<k_3\le 700$ and all $l_i$, $i=1,2,3$, satisfying the constraints. 

We define a \textit{cube} in the $d$-dimensional Euclidean space $\er^d$ to be any translate of the unit cube $[0,1)^d$. Let $T$ be a subset of $\er^d$. The family  
$[0,1)^{d}+T:= \{[0,1)^{d}+t\colon t\in T\}$ is said to be  a \textit{cube tiling} of $\mathbb R^{d}$ if for each pair of distinct vectors $s,t\in T$ the cubes $[0,1)^{d}+s$ and $[0,1)^{d}+t$ are disjoint and
$\bigcup[0,1)^{d}+T=\mathbb R^{d}$. Let $\gruba k :=(k_1,\ldots,k_d)$ be a vector with all coordinates that are positive integers. The tiling $[0,1)^{d}+T$ is said to be $\gruba k$-\textit{periodic} if for every vector of the standard basis $e_1=(1,0,\ldots,0),\ldots, e_d=(0,\ldots,1)$ one has 
$$
T+k_ie_i=T.
$$

We define the \textit{(flat) torus} $\te_{\gruba k}^d$, to be the set 
$[0,k_1)\times\cdots\times[0,k_d)$ with addition $\przyst \gruba k$:
$$
x\oplus y:=((x_1+y_1)\przyst {k_1},\ldots, (x_d+y_d)\przyst{k_d}).
$$
We can extend the notion of a cube so that it will apply to flat tori: Cubes in $\te_{\gruba k}^d$
are the sets of the form $[0,1)^d\oplus t$, where $t\in \te_{\gruba k}^d$. It is clear that we can speak about cube tilings
of $\te_{\gruba k}^d$ and that there is a canonical `one-to-one' correspondence between these tilings and the $\gruba k$-periodic tilings of $\er^d$.   

From now on, $\gruba k$ is assumed to  have all coordinates greater than 1.

If $[0,1)^d\oplus T$ is a cube tiling of $\te_{\gruba k}^d$ and $S\subseteq T$, then we say that the packing $[0,1)^d\oplus S$ is a \textit{simple component} of the cube tiling  if $S$ is an equivalence class of the relation `$\sim$' defined on $T$ as follows:   
$$
\text{$x\sim y$ if and only if $x-y\in \zet^d$.}
$$   

For each $t\in \te_{\gruba k}^d$, the  \textit{integer code} of $t$ is defined by 
$$
\varepsilon(t):=(\lfloor t_1\rfloor +1,\ldots,\lfloor t_d\rfloor+1).
$$ 
Clearly, $\varepsilon$ maps $\te_{\gruba k}^d$ into $[\gruba k]$. 

One can prove the following rather non-trivial result.
\begin{tw}
If $[0,1)^d+S$ is a simple component of a cube tiling $[0,1)^d\oplus T$ of $\te_{\gruba k}^d$, then $\varepsilon (S)\subseteq [\gruba k]$ is complementable by lines.     

If $D\subseteq [\gruba k]$ is complementable by lines, then there is a cube tiling $[0,1)^d\oplus T$ of $\te_{\gruba k}^d$ with a simple component $[0,1)^d+S$  such that $\varepsilon(S)=D$. 
\end{tw}
 
If, in addition, we take into account that by Keller's theorem (see any of the three papers we refer to), the restriction $\varepsilon |T$ of $\varepsilon$ to $T$ is a bijection for each cube tiling $[0,1)^d\oplus T$, then the conjecture can be rephrased as follows: 
\begin{itemize}
\item[] 
If $[0,1)^d+S$ is a simple component of the cube tiling $[0,1)^d+T$ of the torus $\te_{\gruba k}^d$, then the size $|S|$ of $S$ is representable by $\gruba k$. 
\end{itemize}


\begin{thebibliography}{ABCDEF}
{\small
\bibitem {IP} \textsc{A. Iosevich, S. Pedersen}, Spectral and Tiling
Properties of the Unit Cube, \textit{Inter. Math. Res. Notices} {\bf 16} (1998), 819--828.
\bibitem {K} \textsc{O.-H. Keller},  \"Uber die l\"uckenlose Erf\"ullung des Raumes mit
W\"urfeln, \textit{J. Reine Angew. Math.} {\bf 163} (1930), 231--248.
\bibitem {P} \textsc{O. Perron}, \"Uber l\"uckenlose Ausf\"ullung des $n$-dimensionalen Raumes durch
kongruente W\"urfel, \textit{Math. Z.} {\bf 46} (1940), 1--26.
}
\end{thebibliography}
\end{document}